\documentstyle[
11pt,amsfonts]{article}

\newcommand {\rel} {{\mathbb R}}
\newcommand {\com} {{\mathbb C}}
\newcommand {\nat} {{\mathbb N}}

\newtheorem{proposition}{Proposition}[section]
\newtheorem{theorem}{Theorem}[section]

\newcommand{\nc}{\newcommand}

\newcommand{\bcite}[1] {\cite{#1}}

\newcommand {\pr} {\bf}

\newcommand{\proof} {   \begin{flushright}
                        ///
                        \end{flushright}
                }

\newcommand{\defin} { \hspace*{\fill} $\Box$ }

\newcommand {\hyp}  {{\mathcal{H}}}

\def    \mean   {{ \vec {\bf H} }}
\def    \meansc {{ H }}

\def	\d	{{ \ {\rm d} }}

\def    \W      {{ {\mathcal W} }}

\def	\N	{{ \cal N }}

\def	\H	{{ \cal H }}

\nc{\energ}[1]  {{ e_{#1} }}

\def	\poin	{{ poin }}

\def	\geu	{{ g_{euc} }}

\nc{\gpo}[1]	{{ g_{\poin}^{#1} }}
\nc{\gpu}[1]	{{ g_{\poin,{#1}} }}

\nc{\gpmo}[1]	{{ g_{\poin,m}^{#1} }}
\nc{\tgpmo}[1]	{{ \tilde g_{\poin,m}^{#1} }}

\nc{\gplo}[1]	{{ g_{\poin,\lambda}^{#1} }}

\nc{\conf}[1]   {{ [{#1}] }}

\def	\mani	{{ \Gamma }}

\def	\iso	{{ \cal A }}

\nc{\doub}[1]{{ \ddot{#1} }}
\nc{\dd}{ \begin{displaymath} }
\nc{\df}{ \end{displaymath} }
\nc{\dcd}{ \begin{displaymath} \begin{array}{c}}
\nc{\dcf}{ \end{array} \end{displaymath} }
\nc{\ee}{ \begin{equation} }
\nc{\ef}{ \end{equation} }
\nc{\ad}{ \begin{array}{c} }
\nc{\af}{ \end{array} }

\nc{\add}[1]    {}

\begin{document}

\begin{center}
{\huge \bf The umbilic set of Willmore surfaces}
\\ \ \\
Reiner M. Sch\"atzle \\
Fachbereich Mathematik der
Eberhard-Karls-Universit\"at T\"ubingen, \\
Auf der Morgenstelle 10,
D-72076 T\"ubingen, Germany, \\
email: schaetz@everest.mathematik.uni-tuebingen.de \\

\end{center}
\vspace{1cm}

\begin{quote}

{\bf Abstract:}
It is well known that the umbilic points of minimal surfaces
in spaces of constant sectional curvature
consist only of isolated points
unless the surface is totally umbilic on some connected component,
as for example the Hopf form is holomorphic.
In this note, we prove that on Willmore surfaces in codimension one
the umbilic set is locally a one dimensional real-analytic manifold
without boundary or an isolated point.
\end{quote}

\vspace{1cm}



\section{Introduction} \label{intro}

For an immersed surface $\ f : \Sigma
\rightarrow \rel^3\ $ the Willmore functional
is defined by
\begin{displaymath}
        \W(f) = \frac{1}{4} \int \limits_\Sigma |\meansc|^2 \d \mu_g,
\end{displaymath}
where $\ \meansc\ $ denotes the scalar mean curvature of $\ f\ $,
$\ g = f^* g_{euc}\ $ the pull-back metric
and $\ \mu_g\ $ the induced
area measure of $\ f \mbox{ on } \Sigma\ $.

Critical points of the Willmore functional,
called Willmore surfaces or immersions,
satisfy the Euler-Lagrange equation
\begin{equation} \label{intro.euler}
	\Delta_g \meansc + |A^0|_g^2 \meansc = 0,
\end{equation}
see i.e. \bcite{kuw.schae.will} \S2,
where the laplacian along $\ f\ $ is used,
$\ A^0 = A - \frac{1}{2} \mean g\ $
is the tracefree second fundamental form of $\ f\ $.
In particular, smooth Willmore immersions are real-analytic,
see \bcite{morrey.ana}.

A point is called umbilic,
when the tracefree part of the second fundamental form $\ A^0\ $ vanishes,
likewise the set of non-umbilic points is given by
\begin{equation} \label{intro.umbilic}
	\N := [A^0 \neq 0].
\end{equation}
As Willmore immersions are real-analytic,
$\ \N = [A^0 \neq 0] \mbox{ is dense in } \Sigma\ $,
if $\ f\ $ is not totally umbilic on any connected component.

Minimal surfaces are particular examples of Willmore surfaces by (\ref{intro.euler}),
and by conformal invariance of the Willmore functional, see \bcite{chen.conf},
this is also true for minimal surfaces in the sphere $\ S^3\ $
or in the hyperbolic space $\ \hyp^3\ $
after applying a local conformal diffeomorphism to the euclidean space.
We remark that $\ \rel^3, S^3 \mbox{ and } \hyp^3\ $ all
have constant sectional curvature.
Now the umbilic points of minimal surfaces consist only of isolated points
unless the surface is totally umbilic on some connected component.
In codimension one, this can easily be seen by considering
the scalar second fundamental form $\ h^0_{ij} := \langle A^0_{ij} , \nu \rangle
\mbox{ for a smooth unit normal } \nu \mbox{ of } f\ $,
the function
\begin{equation} \label{intro.hopf}
	\varphi := h^0_{11} - i h^0_{12}
\end{equation}
and the quadratic Hopf form $\ \H = (\varphi/2) (dz)^2\ $.
The Hopf form is defined independent of the oriented conformal local chart
and changes to its conjugate when switching the orientation of the local chart.
In case of minimal surfaces, $\ \varphi\ $ is holomorphic,
see for example \bcite{lawson.mini} Lemma 1.2 for the sphere,
hence its zeros, which are precisely the umbilic points,
consist only of isolated points
unless $\ \varphi\ $ vanishes identically on some connected component.

The aim of this article is to show
that on Willmore surfaces in codimension one
the umbilic set is locally a one dimensional real-analytic manifold
without boundary or an isolated point.
\\ \ \\
{\bf Theorem \ref{struc.theo}}
{\it
For any smooth not totally umbilic Willmore immersion
$\ f: \Sigma \rightarrow \rel^3\ $
of an open connected surface $\ \Sigma\ $,
the set of umbilic points is a closed set in $\ \Sigma\ $
which is locally a one dimensional real-analytic manifold
of $\ \Sigma\ $ without boundary
or an isolated point.

Therefore the set of umbilic points of $\ f\ $ can be written
\begin{displaymath}
	\Sigma - \N = [A^0 = 0]
	= \mani + \iso,
\end{displaymath}
where $\ \mani\ $ is a closed set
and a one dimensional real-analytic submanifold
of $\ \Sigma\ $ without boundary
and $\ \iso\ $ is a closed set isolated points.
}
\defin


\setcounter{equation}{0}

\section{Cauchy-Riemann equation} \label{cr}

We consider a smooth immersion $\ f: \Sigma \rightarrow \rel^3
\mbox{ of a surface } \Sigma\ $,
or more precisely in a local conformal chart
a conformal immersion $\ f: B_1(0) \subseteq \com \rightarrow \rel^3\ $
with pull-back metric $\ g = f^* \geu = e^{2u} \geu\ $
and use the differential operators of the Wirtinger calculus
\begin{displaymath}
        \partial_z = \frac{1}{2} (\partial_1 - i \partial_2)
        \quad \mbox{and} \quad
        \partial_{\bar z} = \frac{1}{2} (\partial_1 + i \partial_2).
\end{displaymath}
We derive a Cauchy-Riemann equation for the pair of functions
$\ \varphi\ $, defined in (\ref{intro.hopf}), and $\ \partial_z \meansc\ $,
and we start with an equation for $\ \varphi\ $.

\begin{proposition} \label{cr.hopf}

For any conformal immersion $\ f: B_1(0) \subseteq \com \rightarrow M,
M \in \{ \rel^3, S^3, \hyp^3 \}\ $,
with pull-back metric $\ g = f^* g_M = e^{2u} \geu\ $,
we have
\begin{equation} \label{cr.hopf.deri}
	\partial_{\bar z} \varphi = \frac{e^{2u}}{2} \partial_z \meansc,
\end{equation}
where $\ \meansc\ $ is the scalar mean curvature of $\ f \mbox{ in } M\ $.
\end{proposition}
{\pr Beweis:} \\
From (\ref{intro.hopf}), we calculate
\begin{displaymath}
	2 \partial_{\bar z} \varphi
	= (\partial_1 h^0_{11} + \partial_2 h^0_{12})
	- i (\partial_1 h^0_{12} - \partial_2 h^0_{11})
\end{displaymath}
\begin{displaymath}
	= (\partial_1 h^0_{11} + \partial_2 h^0_{21})
	- i (\partial_1 h^0_{12} + \partial_2 h^0_{22}) =
\end{displaymath}
\begin{equation} \label{cr.hopf.aux}
	= e^{2u} g^{kl} (\partial_k h^0_{l1} - i \partial_k h^0_{l2}),
\end{equation}
as $\ h^0\ $ is symmetric and tracefree with respect to $\ g = e^{2u} \geu\ $.
The covariant derivatives with respect to $\ g = e^{2u} \geu\ $ are defined by
\begin{displaymath}
	\nabla_k h^0_{lm}
	= \partial_k h^0_{lm}
	- \Gamma_{kl}^r h^0_{rm}
	- \Gamma_{km}^r h^0_{lr},
\end{displaymath}
where $\ \Gamma\ $ denote the Christoffel symbols with respect to the metric $\ g = e^{2u} \geu\ $
and are given by
\begin{displaymath}
	\Gamma_{kl}^r
	= \frac{1}{2} g^{rs} (\partial_k g_{ls} + \partial_l g_{sk}
	- \partial_s g_{kl}) =
\end{displaymath}
\begin{displaymath}
	= g^{rs} (g_{ls} \partial_k u + g_{sk} \partial_l u
	- g_{kl} \partial_s u)
	= \delta_l^r \partial_k u
	+ \delta_k^r \partial_l u - g_{kl} g^{rs} \partial_s u
\end{displaymath}
for the Kronecker symbols $\ \delta\ $.
We calculate
\begin{displaymath}
	g^{kl} \Gamma_{kl}^r
	= g^{kl} (\delta_l^r \partial_k u
	+ \delta_k^r \partial_l u - g_{kl} g^{rs} \partial_s u) =
\end{displaymath}
\begin{equation} \label{cr.hopf.gamma1}
	= g^{kr} \partial_k u
	+ g^{lr} \partial_l u - 2 g^{rs} \partial_s u = 0
\end{equation}
and
\begin{displaymath}
	g^{kl} \Gamma_{km}^r h^0_{lr}
	= g^{kl} (\delta_m^r \partial_k u
	+ \delta_k^r \partial_m u - g_{km} g^{rs} \partial_s u) h^0_{lr} =
\end{displaymath}
\begin{displaymath}
	= g^{kl} h^0_{lm} \partial_k u
	+ g^{kl} h^0_{lk} \partial_m u - \delta_m^l g^{rs} h^0_{lr} \partial_s u =
\end{displaymath}
\begin{equation} \label{cr.hopf.gamma2}
	= g^{kl} h^0_{lm} \partial_k u
	- g^{rs} h^0_{rm} \partial_s u = 0,
\end{equation}
as $\ h^0\ $ is symmetric and tracefree with respect to $\ g\ $.
Together
\begin{displaymath}
	g^{kl} \nabla_k h^0_{lm}
	= g^{kl} \partial_k h^0_{lm}
\end{displaymath}
and plugging into (\ref{cr.hopf.aux})
\begin{displaymath}
	2 \partial_{\bar z} \varphi
	= e^{2u} g^{kl} (\nabla_k h^0_{l1}
	- i \nabla_k h^0_{l2}).
\end{displaymath}
On the other hand by the Mainardi-Codazzi equation,
as $\ \rel^3, S^3, \hyp^3\ $ all have constant sectional curvature,
see \bcite{docarmo} \S 6 Proposition 3.4 and \S 4 Lemma 3.4,
\begin{displaymath}
	g^{kl} \nabla_k h^0_{lm}
	= g^{kl} \nabla_k (h_{lm} - \frac{1}{2} g_{lm} \meansc)
	= g^{kl} (\nabla_m h_{kl} - \frac{1}{2} g_{lm} \partial_k \meansc)
	= \frac{1}{2} \partial_m \meansc,
\end{displaymath}
hence
\begin{displaymath}
	\partial_{\bar z} \varphi
	= \frac{e^{2u}}{4} (\partial_1 \meansc - i \partial_2 \meansc) 
	= \frac{e^{2u}}{2} \partial_z \meansc,
\end{displaymath}
which is (\ref{cr.hopf.deri}).
\proof
{\large \bf Remark:} \\
For minimal surfaces in the constant curvature spaces $\ \rel^3, S^3 \mbox{ or } \hyp^3\ $,
we get from (\ref{cr.hopf.deri}) that $\ \varphi\ $ is holomorphic,
hence the umbilic points, which are the zeros of $\ \varphi\ $,
consist only of isolated points,
if $\ f\ $ is not totally umbilic on any connected component.
\defin

\begin{proposition} \label{cr.cr}

For any conformal Willmore immersion $\ f: B_1(0) \subseteq \com \rightarrow \rel^3\ $,
we have the Cauchy-Riemann equation
\begin{equation} \label{cr.cr.equ}
	\partial_{\bar z} (\varphi , \partial_z\meansc)^T
	= M \cdot (\varphi , \partial_z \meansc)^T
\end{equation}
for some smooth matrix $\ M: B_1(0) \rightarrow \com^{2 \times 2}\ $.
\end{proposition}
{\pr Proof:} \\
(\ref{cr.hopf.deri}) implies the first row of (\ref{cr.cr.equ}).
Now $\ f\ $ satisfies as Willmore immersion
the Euler-Lagrange equation (\ref{intro.euler}),
hence by the definition in (\ref{intro.hopf}) that
\begin{displaymath}
	\Delta_g \meansc	
	= -|A^0|_g^2 \meansc
	= -g^{jk} g^{lm} \langle A^0_{jl} , A^0_{km} \rangle \meansc =
\end{displaymath}
\begin{displaymath}
	= -e^{-4u} \Big( (h^0_{11})^2 + (h^0_{12})^2 + (h^0_{21})^2 + (h^0_{22})^2 \Big) \meansc
	= -2 e^{-4u} |\varphi|^2 \meansc,
\end{displaymath}
as $\ h^0\ $ is symmetric and tracefree with respect to $\ g\ $.
On the other hand using (\ref{cr.hopf.gamma1}), we get
\begin{displaymath}
	\Delta_g \meansc
	= g^{kl} \nabla_k \nabla_l \meansc
	= g^{kl} \partial_k \partial_l \meansc
	- g^{kl} \Gamma_{kl}^m \partial_m \meansc =
\end{displaymath}
\begin{displaymath}
	= e^{-2u} (\partial_1 \partial_1 \meansc
	+ \partial_2 \partial_2 \meansc)
	= 4 e^{-2u} \partial_{\bar z} \partial_z \meansc.
\end{displaymath}
Together
\begin{displaymath}
	\partial_{\bar z} \partial_z \meansc
	= \frac{1}{4} e^{2u} \Delta_g \meansc
	= -\frac{1}{2} e^{-2u} \bar \varphi \meansc \varphi,
\end{displaymath}
which gives the second row of (\ref{cr.cr.equ}).
\proof


\setcounter{equation}{0}

\section{The structure of the umbilic set} \label{struc}

In this section we prove our main theorem.

\begin{theorem} \label{struc.theo}

For any smooth not totally umbilic Willmore immersion
$\ f: \Sigma \rightarrow \rel^3\ $
of an open connected surface $\ \Sigma\ $,
the set of umbilic points is a closed set in $\ \Sigma\ $
which is locally a one dimensional real-analytic manifold
of $\ \Sigma\ $ without boundary
or an isolated point.

Therefore the set of umbilic points of $\ f\ $ can be written
\begin{equation} \label{intro.struc.decomp}
	\Sigma - \N = [A^0 = 0]
	= \mani + \iso,
\end{equation}
where $\ \mani\ $ is a closed set
and a one dimensional real-analytic submanifold
of $\ \Sigma\ $ without boundary
and $\ \iso\ $ is a closed set isolated points.
\end{theorem}
{\pr Proof:} \\
By continuity, the set of umbilic points $\ [A^0 = 0] \mbox{ is closed in } \Sigma\ $.
We consider a conformal Willmore immersion
$\ f: B_1(0) \subseteq \com \rightarrow \rel^3\ $
with pull-back metric $\ g = f^* \geu = e^{2u} \geu\ $
which has $\ 0\ $ as an umbilic point.
Clearly $\ [A^0 = 0] = [\varphi = 0]\ $.

By the Cauchy-Riemann equation for $\ (\varphi , \partial_z \meansc)\ $ in Proposition \ref{cr.cr}
and observing that $\ \varphi\ $ is real-analytic
and does not vanish identically,
as $\ f\ $ is real-analytic and not totally umbilic,
there exists by \bcite{esch.tri.branch} Lemmas 2.1 and 2.2
an integer $\ m \in \nat_0\ $ such that
$\ (\varphi,\partial_z \meansc)(z) = z^m (\psi,\chi)(z)\ $
for some smooth $\ \psi, \chi: B_1(0) \rightarrow \com
\mbox{ with } (\psi(0),\chi(0)) \neq 0\ $.
As $\ \varphi , \partial_z \meansc\ $ are real-analytic,
$\ \psi, \chi\ $ are real-analytic up to the origin as well.

As the origin is considered to be an umbilic point, we have
\begin{displaymath}
	[A^0 = 0] = [\varphi = 0] = [\psi = 0] \cup \{0\}.
\end{displaymath}
If $\ \psi(0) \neq 0\ $,
the origin is an isolated umbilic point.

If $\ \psi(0) = 0\ $, then $\ \chi(0) \neq 0\ $.
Observing by (\ref{cr.hopf.deri}) that
\begin{displaymath}
	z^m \partial_{\bar z} \psi(z) = \partial_{\bar z} \varphi(z)
	= \frac{e^{2u(z)}}{2} \partial_z \meansc(z)
	= z^m \frac{e^{2u(z)}}{2} \chi(z),
\end{displaymath}
we get by continuity $\ \partial_{\bar z} \psi = (e^{2u}/2) \chi\ $ and
\begin{displaymath}
	\partial_{\bar z} \psi(0) = \frac{e^{2u(0)}}{2} \chi(0) \neq 0.
\end{displaymath}
Therefore $\ D \psi(0) \neq 0\ $
or likewise $\ \nabla Re(\psi)(0) \neq 0
\mbox{ or } \nabla Im(\psi)(0) \neq 0\ $.
Then by the implicit function theorem,
the set $\ [Re(\psi) = 0] \mbox{ or the set } [Im(\psi) = 0]\ $
is a real-analytic curve locally around the origin,
and the set of umbilic points
\begin{displaymath}
	[\varphi = 0] = [\psi = 0]
	= \Big( [Re(\psi) = 0] \cap [Im(\psi) = 0] \Big)
\end{displaymath}
is contained in this real-analytic curve locally around the origin.
Since $\ \varphi\ $ is real-analytic,
we see that either the origin is an isolated umbilic point
or the whole real-analytic curve belongs in a neighbourhood
to the set of umbilic points,
and the theorem is proved.
\proof



\end{document}